\documentclass[12pt]{amsart}

\usepackage{amsfonts, amssymb, amscd}
\usepackage{verbatim}
\usepackage{eucal}
\usepackage{amssymb}
\usepackage{mathrsfs}
\usepackage{graphicx}
\usepackage{psfrag}

\def\ZZ{{\mathbb Z}}

\def\PP{{\mathbb P}}
\def\CC{{\mathbb C}}

\addtocounter{footnote}{1}

\newtheorem{lemma}{Lemma}[section]

\newtheorem{corollary}[lemma]{Corollary}

\newtheorem{proposition}[lemma]{Proposition}
\theoremstyle{definition}

\theoremstyle{remark}
\newtheorem*{proof*}{Proof}
\numberwithin{equation}{section}

\title[On Calabi-Yau threefolds]{On Calabi-Yau threefolds with 
large nonabelian fundamental groups}

\author{Lev Borisov and Zheng Hua}
\address{Department of Mathematics \\ University of Wisconsin \\
  Madison \\ WI \\ 53706 \\ USA\\{\tt borisov@math.wisc.edu}\\
{\tt hua@math.wisc.edu}  }

\thanks{Lev Borisov has been partially supported
by the National Science Foundation under grant No.\ DMS-0456801.}

\begin{document}

\begin{abstract}
In this short note we construct Calabi-Yau threefolds with
nonabelian fundamental groups of order $64$
as quotients of the small resolutions
of certain complete intersections of quadrics in $\PP^7$ that were 
first considered by M. Gross and S. Popescu.
\end{abstract}

\maketitle

\section{Introduction}
In \cite{GPo}, M. Gross and S. Popescu constructed 
Calabi-Yau varieties that admit a free action of the
group $(\ZZ/8\ZZ)^2$. We observe that these varieties 
admit a larger automorphism group, and that two other
subgroups of order $64$ of this group act freely.
Our construction is very simple, modulo the nontrivial
calculations of \cite{GPo}. The essential idea is
that in the presence of only two, nonisomorphic, minimal models
every automorphism of the singular variety $V_{8,y}$
naturally extends to the models. We then construct
the additional automorphisms explicitly as permutations
of coordinates.

In the last section we very briefly describe what led to the
discovery of these examples and how they fit into 
a more general question of finding free group actions
on complete intersections of quadrics.

{\bf Acknowledgements.} We have used extensively
GAP software package in our project, although no computer
calculations are necessary to check the results
of this paper.

\section{Construction}
We recall the definition of the Calabi-Yau varieties
of Gross and Popescu, \cite{GPo,GPa}. 
Consider the projective space 
$\CC\PP^7$ with coordinates $(x_0:\ldots:x_7)$.
For any $(y_0:y_1:y_2)\in \CC\PP^2$ in general position,
consider the intersection $V_{8,y}$  of four quadrics
$$
\begin{array}{l}
y_1y_3(x_0^2+x_4^2) - y_2^2(x_1x_7+x_3x_5)+(y_1^2+y_3^2)x_2x_6=0\\
y_1y_3(x_1^2+x_5^2) - y_2^2(x_2x_0+x_4x_6)+(y_1^2+y_3^2)x_3x_7=0\\
y_1y_3(x_2^2+x_6^2) - y_2^2(x_3x_1+x_5x_7)+(y_1^2+y_3^2)x_4x_0=0\\
y_1y_3(x_3^2+x_7^2) - y_2^2(x_4x_2+x_6x_0)+(y_1^2+y_3^2)x_5x_1=0.\\
\end{array}
$$
This variety has $64$ ODP singular points which are the orbit of 
$(0:y_1:y_2:y_3:0:-y_3:-y_2:-y_1)$ under the action of the group
$G=(\ZZ/8\ZZ)^2$ generated by $\tau$ and $\sigma$ where $\tau(x_i)=
\zeta_8^{-i}x_i$ and $\sigma$ is a cyclic permutation of the variables 
given by the cycle
$(01234567)$. This group $G$ acts freely on $V_{8,y}$.

The variety $V_{8,y}$ admits two small resolutions $V_{8,y}^1$
and $V_{8,y}^2$ both of which are Calabi-Yau threefolds with 
$h^{1,1}=h^{1,2}=2$. These resolutions are related by a flop and 
are obtained from $X$ by blowups of abelian surfaces in $X$ of 
degrees $32$ and $16$ respectively. 

The key to the current note is the following lemma that relies
heavily on the  results of \cite{GPo}.
\begin{lemma}\label{2.1}
Every automorphism $\gamma$ of $V_{8,y}$ lifts
to automorphisms of $V_{8,y}^1$ and $V_{8,y}^2$. 
If $\gamma$ acts freely on $V_{8,y}$, then its lifts act freely
on the small resolutions.
\end{lemma}

\begin{proof}
We will show that $\gamma$ lifts to $V_{8,y}^2$, with the other 
case being completely analogous.

As in \cite{GPo}, we can think of $V_{8,y}^2$ as the blowup of $V_{8,y}$
by the ideal of some Weil divisor $A$. Then there is a natural
isomorphism $\hat\gamma$ from $V_{8,y}^2$ to the blowup of $V_{8,y}$ by the ideal of $\gamma(A)$. Denote the latter blowup by $Z$. It is 
clearly a minimal model of $V_{8,y}$, so by the observation of 
\cite{GPo}, $Z$ is isomorphic to either $V_{8,y}^2$ or $V_{8,y}^1$
(over $V_{8,y})$. To prove the first assertion of the lemma, it 
remains to show that  $V_{8,y}^2$ and $V_{8,y}^1$ are not 
isomorphic, since we could then compose $\hat\gamma$
with an isomorphism $Z\to V_{8,y}^2$ over $V_{8,y}$. 
From the description of their K\"ahler cones 
we can see that any such isomorphism between
$V_{8,y}^{i}$ would have to respect 
their structures as fibrations over $\PP^1$ as well as their maps to
$\PP^7$. However, the degrees of the fibers are different in these 
two cases, hence $V_{8,y}^2\not\cong V_{8,y}^1$.

To show the last assertion, observe that if the lift of $\gamma$ had a 
fixed point, then the image of that point in $V_{8,y}$ would be
fixed by $\gamma$.
\end{proof}

We can now describe our construction. 
\begin{proposition}\label{2.2}
We define the group $G_1$ generated by $\tau$ and 
the permutation $\sigma_1=(07214365)$ of the coordinates $x_i$.
Then $G_1$ is a nonabelian group isomorphic 
to a semidirect product of two copies of $\ZZ/8\ZZ$.
We define the group $G_2$ generated by $\tau$ and 
the permutations
$\sigma_2= (0246)(1357)$ and $\sigma_3=(0145)(3276)$.
Then $G_2$ is a nonabelian group isomorphic 
to a semidirect product of normal subgroup
$\ZZ/8\ZZ$ generated by $\tau$ and the quaternion group
$H$ generated by $\sigma_2$ and $\sigma_3$.
Both $G_1$ and $G_2$ act freely on $V_{8,y}$.
\end{proposition}

\begin{proof}
The structure of $G_1$ and $G_2$ is immediate from
their definition, and 
it is straightforward to see that $\sigma_i$ acts on $V_{8,y}$.
To show that they act freely, it is enough to check the action of 
all involutions in $G_i$. In view of 
$1\to \langle \tau\rangle\to G_1\to \ZZ/8\ZZ\to 1$,
all involutions of $G_1$ lie in the subgroup generated by 
$\tau$ and $\sigma_1^4$.  Similarly, all involutions of 
$G_2$ lie in the subgroup generated by $\tau$ and $\sigma_2^2$.
It remains to observe that 
these subgroups are contained in $G$, which is known to act freely \cite{GPo}.
\end{proof}

\begin{corollary}
For $i=1,2$, $j=1,2$, the quotients of $V_{8,y}^i$ by
(the lift of) the group $G_j$ are smooth Calabi-Yau threefolds
with fundamental groups $G_j$ of order $64$.
\end{corollary}

\begin{proof} By \cite{GPa}, $V_{8,y}^{i}$ are simply connected.
It remains to combine Lemma \ref{2.1} and Proposition \ref{2.2}.
\end{proof}

We remark that $G_2$ contains the quaternion group in its
regular representation, and $V_{8,y}$ can be thought of
as a singular member of the family constructed by Beauville
in \cite{B}. Also, the quotients $V^i_{8,y}/G_{j}$ have
Hodge numbers $h^{1,1}=h^{1,2}=2$ and the structures
of fibration with abelian surface fibers.

\section{Comments}
We have stumbled upon these examples largely by chance.
We originally set out to investigate free actions of finite
groups on complete intersections of four quadrics in hopes
of extending the construction of \cite{B}. Every such group 
action naturally leads to a projective representation 
of dimension eight, which can be then thought of as 
a linear representation of a Schur cover $S$ of the finite group.
If the action is free,  the holomorphic Lefschetz formula
essentially dictates what character of $S$ one needs to consider,
moreover, it leads to strong restrictions on possible group
actions. We have used the GAP software package extensively in
our search.

It turns out that the maximum possible order of the 
group (if one allows the variety to have ODP 
singularities) is $64$. We have found five possible groups of order $64$.
Upon closer consideration, it turned out that three of them act
on the same family of varieties and they are precisely 
the groups $G$, $G_1$ and $G_2$ that 
appear in this note. The other two, together with the classification
of groups of smaller order will be addressed by the second author
in \cite{H}.

\end{document}